\documentclass{article}
\usepackage{amssymb}
\begin{document}

\begin{center}
{\small Journal of Algebra 252 (2002) 65-73}
\bigskip

{\large\bf Tensor products of Cohen-Macaulay rings}\\
{\bf Solution to a problem of Grothendieck} \bigskip\bigskip

S. BOUCHIBA\\
Department of Mathematics \\
University Moulay Ismail \\
Meknes 50000, Morocco\\
Email address: bouchiba@fsmek.ac.ma  \bigskip\bigskip

\renewcommand{\thefootnote}{\fnsymbol{footnote}}
S. KABBAJ \\
Department of Mathematics \\
Harvard University \\
Cambridge, MA 02138, USA\\
Email address: kabbaj@math.harvard.edu
\end{center}    \vspace{1cm}

\noindent {\small \bf ABSTRACT}

In this paper we solve a problem, originally raised by Grothendieck, on the
transfer of Cohen-Macaulayness to tensor products of algebras over a field $k$. As a prelude to this, we investigate the grade for some specific types of ideals that play a primordial role within the ideal structure of
such constructions.
\vspace{2cm}

\begin{center}INTRODUCTION \end{center}

All rings and algebras considered in this paper are commutative with
identity elements and, unless otherwise specified, are to be assumed to be
non-trivial. All ring homomorphisms are unital. Throughout, $k$ stands for a field. Let $A$ be a ring. We shall use $G(I)$ to denote the grade of an ideal $I$
of $A$, $Z(A)$ to
denote the set of all zero-divisors of $A$, and $k_A(p)$ to denote the
quotient field of
$\frac {A}{p}$ for any prime ideal $p$ of $A$.

Let $A$ be a Noetherian ring and $I$ a proper ideal of $A$. The grade of $I$ is defined to be the common length of all maximal $A$-sequences in $I$. It can be measured by the\\ (non-) vanishing of certain Ext modules. In fact, according to \cite[Theorem 16.7]{m}, $G(I)=$ Inf$\{i|$Ext$_A^i(\frac AI,A)\neq 0\}$. This connection opened commutative algebra to the
application of homological methods. Finally, recall that the Cohen-Macaulay rings are those Noetherian rings in which grade and height coincide for every ideal.

Our aim in this paper is to prove that the Cohen-Macaulay property
is inherited by tensor products of $k$-algebras.
To this purpose, the first section investigates the grade of three specific
types of ideals that play a primordial role within the ideal structure of the tensor product of two $k$-algebras.
This allows us, in the second section, to establish the main theorem, that
is,  for $k$-algebras $A$
and $B$ such that $A\otimes _kB$ is Noetherian, $A\otimes _kB$ is a
Cohen-Macaulay ring if and only if so are $A$ and $B$.

Suitable background on depth of modules and Cohen-Macaulay
rings is \cite{bh}, \cite{gr}, \cite{k}, and \cite{m}. For a geometric treatment of the Cohen-Macaulay property,
we refer the reader to the excellent book of Eisenbud \cite{e}.
Recent developments on heights of primes and dimension theory in tensor products of $k$-algebras are to be found in \cite{bgk1}, \cite{bgk2}, and \cite{w}. Any unreferenced material is standard, as in \cite{g}, and \cite{k}.\bigskip\bigskip

\begin{center} 1. GRADE OF IDEALS IN A TENSOR PRODUCT OF TWO $k$-ALGEBRAS
\end{center}

The grade of an arbitrary ideal in a (Noetherian) tensor product of two
$k$-algebras seems to be difficult to grasp. It would appeal to new techniques yet to be discovered. Our goal here is much more modest. We shall determine the grade of three specific types of ideals that play a primordial role within the ideal structure of this construction.\bigskip

We announce the main result of this section.\bigskip

\noindent {\bf Theorem 1.1.} Let $A$ and $B$ be $k$-algebras such that $
A\otimes _kB$ is Noetherian. Let $I$ and $J$ be proper ideals of $A$ and
$B$,
respectively. Then:\bigskip

{\bf a)} $G(I\otimes _kB)=G(I)$ and, similarly, $G(A\otimes _kJ)=G(J)$.\bigskip

{\bf b)} $G(I\otimes _kB+A\otimes _kJ)=G(I)+G(J)$.\bigskip

{\bf c)} $G(I\otimes _kJ)=$ Inf$(G(I),G(J))$.\bigskip

Let $A$ and $B$ be two $k$-algebras. Let $x$ be a non zero-divisor element
of $A$ and $y$ a non zero-divisor
element of $B$. Then $x\otimes y$ is a non zero-divisor element of $A\otimes
_kB$.
Let $I$ be a proper ideal of $A$. Then $I\otimes _kB$ is a proper ideal of
$A\otimes _kB$. If $x_1,...,x_n$ is an $A$-sequence, then it is easily seen that $x_1\otimes1,...,x_n\otimes1$ is an $(A\otimes_kB)$-sequence.
These basic facts will be used frequently in the sequel without explicit
mention.
Moreover, we assume familiarity with the natural isomorphisms for tensor
products,
as in  \cite{b}. In particular, we identify $A$ and $B$ with their
respective
images in $A\otimes_kB$, and if $I$ and $J$ are proper ideals of $A$ and $B$, respectively, then $\frac {A\otimes_kB}{I\otimes_kB+A\otimes_kJ}\cong \frac AI \otimes_k\frac BJ$. Also, we recall that
$A\otimes_kB$ is a free (hence flat) extension of $A$ and $B$.
\bigskip

The proof of the main theorem requires the following preparatory lemma.\bigskip

Recall first that if $A$ is a ring and $x_1,...,x_n$ are elements of $A$, then $x_1,...,x_n$ is said to be a permutable $A$-sequence if any permutation of the $x$'s is also an $A$-sequence.\bigskip

\noindent {\bf Lemma 1.2.} Let $A$ and $B$ be two $k$-algebras. Let $x_1,...,x_n$ be a permutable $A$-sequence and $y_1,...,y_n$ be a
permutable $B$-sequence. Then $x_1\otimes y_1,...,x_n\otimes y_n$
is a permutable $(A\otimes_kB)$-sequence.\bigskip

\noindent {\it Proof.} The argument follows easily from the
combination
 of the next two statements. The proofs of these are straightforward and hence
 left to the reader.\bigskip

1) If $x_1,x_2,...,x_n$ is a permutable $A$-sequence, then so is $x_1x_2,x_3,...,x_n$.\bigskip

2) If $x_1,...,x_n$ is a permutable $A$-sequence and $y_1,...,y_m$ is a permutable $B$-sequence, then $x_1\otimes 1,...,x_n\otimes 1,1\otimes y_1,...,1\otimes y_m$ is a permutable $(A\otimes_kB)$-sequence. $\Box
$\bigskip

\noindent {\it Proof of the theorem.} a) Let $G(I)=n$ and $x_1,...,x_n$ be
an $A$-sequence in $I$. Then $x_1,...,x_n$ is an $A\otimes
_kB$-sequence in $I\otimes _kB$. Since $I\subseteq Z(\frac
A{(x_1,...,x_n)})$,
there exists $a\in A \setminus (x_1,...,x_n)$ such that $Ia\subseteq
(x_1,...,x_n)$ \cite[Theorem 82]{k}. Then
\[ \begin{array}{ll}
(I\otimes _kB)a  & = Ia\otimes_kB \\
                 & \subseteq (x_1,...,x_n)\otimes _kB \\
& =(x_1,...,x_n).
\end{array} \]
Clearly, $a\notin
(x_1,...,x_n)\otimes _kB$. Hence $I\otimes _kB\subseteq Z(\frac{A\otimes _kB
}{(x_1,...,x_n)})$. Consequently, $G(I\otimes _k~B)=n=G(I)$.
Likewise for $G(A\otimes _kJ)=G(J)$.\bigskip

b) Let $G(I)=n$ and $G(J)=m$. Let $x_1,...,x_n$ be an $A$-sequence in $I$
and $y_1,...,y_m$ a $B$-sequence in $J$. Obviously, $x_1,...,x_n,y_1,...,y_m$ is an
$A\otimes_kB$-sequence in $I\otimes _kB+A\otimes _kJ$. Since $I\subseteq
Z(\frac A{(x_1,...,x_n)})$ and $J\subseteq Z(\frac B{(y_1,...,y_m)})$, there
exist $a\in A \setminus (x_1,...,x_n)$ and $b\in B \setminus (y_1,...,y_m)$
such that $Ia\subseteq (x_1,...,x_n)$ and $Jb\subseteq (y_1,...,y_m)$. Then
\[ \begin{array}{ll}
(I\otimes _kB+A\otimes _kJ)(a\otimes b) & \subseteq Ia\otimes
_kB+A\otimes_kJb\\
                                        & \subseteq (x_1,...,x_n)\otimes
_kB+A\otimes_k(y_1,...,y_m)\\
                                        & = (x_1,...,x_n,y_1,...,y_m).
\end{array} \]
Since $\overline{a}\neq \overline{0}$ in $\frac
A{(x_1,...,x_n)} $ and $\overline{b}\neq \overline{0}$ in $\frac
B{(y_1,...,y_m)}$, then $\overline{a}\otimes \overline{b}\neq \overline{0}$
in $\frac{A\otimes _kB}{(x_1,...,x_n,y_1,...,y_m)}$, whence
$(I\otimes _kB+A\otimes_kJ)\subseteq Z(\frac{A\otimes
_kB}{(x_1,...,x_n,y_1,...,y_m)})$.
Consequently, $G(I\otimes _kB+A\otimes _kJ)=G(I)+G(J)$, as asserted.\bigskip

c) Let $G(I)=n\leq G(J)=m$. By \cite[Exercise 23, p. 104]{k}, there exist a permutable
$A$-sequence $x_1,...,x_n$ in $I$ and a permutable
$B$-sequence $y_1,...,y_m$ in $J$. By Lemma 1.2,
$x_1\otimes y_1,...,x_n\otimes y_n$ is an $A\otimes_kB$-sequence in
$I\otimes _kJ$.
Since, by (a), $n=G(I\otimes_kB)\geq G(I\otimes_kJ)$, it follows that $G(I\otimes_kJ)=n$, as desired.
$\Box $\bigskip\bigskip

\begin{center} 2. WHEN IS THE TENSOR PRODUCT OF TWO $k$-ALGEBRAS A
COHEN-MACAULAY RING? \end{center}

Recall that a Cohen-Macaulay ring is a Noetherian ring $A$ in which
$G(M)=ht(M)$ for every maximal ideal $M$ of $A$ \cite[Definition, p. 95]{k}. It is worthwhile noting, according to [9, Theorem 136], that grade and height coincide for every proper ideal in a Cohen-Macaulay ring. In
general, the inequality height $\geq$ grade holds. We next show that, in $A \otimes_kB$, the assumption of equality of grade and height for ideals of the form $p\otimes_kB + A\otimes_kq$, where $p$ and $q$ are prime ideals of $A$ and $B$, respectively, implies equality for all ideals.

In 1965, Grothendieck proved in \cite[(6.7.1.1)]{gr} that if $K$ and $L$ are
extension fields of $k$ one of which is finitely generated over $k$, then
$K\otimes _kL$ is a Cohen-Macaulay ring. In 1969, Watanabe et al. extended
this result showing, in \cite[Theorem]{wa}, that if $A$ and $B$
are Cohen-Macaulay rings such that $A\otimes _kB$ is Noetherian and
$\frac Am$ is a finitely generated field extension of $k$ for each maximal ideal $m$ of $A$, then $A\otimes _kB$ is a Cohen-Macaulay ring.

Our purpose in this section is to prove the following:\bigskip

\noindent {\bf Theorem 2.1.} Let $A$ and $B$ be $k$-algebras such that $A
\otimes_kB$ is
Noetherian. Then the following statements are equivalent:\bigskip

{\bf i)} $A \otimes_kB$ is a Cohen-Macaulay ring ;\bigskip

{\bf ii)} $G(I \otimes_kB + A \otimes_kJ) = ht(I \otimes_kB +
A\otimes_kJ)$,
for all proper ideals $I$ and $J$ of $A$ and $B$, respectively ;\bigskip

{\bf iii)} $G(p \otimes_kB + A \otimes_kq) = ht(p \otimes_kB + A
\otimes_kq)$,
for all prime ideals $p$ and $q$ of $A$ and $B$, respectively ;\bigskip

{\bf iv)} $A$ and $B$ are Cohen-Macaulay rings.\bigskip

The discussion which follows, concerning basic facts about $k$-algebras,
will provide
some background to the main theorem and will be of use in its proof. We shall use $t.d.(A:k)$ to denote the transcendence degree of a $k$-algebra $A$ over $k$. It is worth reminding the reader that for an arbitrary $k$-algebra $A$ (not necessarily a domain), $t.d.(A:k):=$ Sup$\{t.d.(\frac Ap:k)|p\in$ Spec$(A)\}$ (cf. \cite[p. 392]{w}).

Notice first that the tensor product of two extension fields of $k$ is not
necessarily Noetherian \cite{v}.
However, given two $k$-algebras $A$ and $B$ such that $A\otimes _kB$ is
Noetherian, then $A$ and $B$ are necessarily Noetherian rings; moreover,
either $t.d.(A:k)<\infty$ or $t.d.(B:k)<\infty$: indeed,
since $A$ and $B$ each have only a finite number of minimal prime ideals, there exist $p\in $ Spec$(A)$ and $q\in $ Spec$(B)$ such that $t.d.(A:k)=t.d.(\frac Ap:k)$ and $t.d.(B:k)=t.d.(\frac Bq:k)$.
Clearly, $k_A(p)\otimes _kk_B(q)$ is Noetherian, since it is a localization of $\frac Ap\otimes _k\frac
Bq\cong \frac{A\otimes _kB}{p\otimes _kB+A\otimes _kq}$, which is
Noetherian. We obtain, by \cite[Corollary 10]{v}, that either $t.d.(k_A(p):k)<\infty $ or
$t.d.(k_B(q):k)<\infty $, as desired.\bigskip

The proof of the main theorem requires two preparatory results.\bigskip

\noindent {\bf Lemma 2.2.} Let $K$ and $L$ be two extension fields of $
k$ such that $K\otimes _kL$ is Noetherian. Then $K\otimes _kL$ is
a Cohen-Macaulay ring.\bigskip

\noindent {\it Proof.} Step 1. We claim that $K\otimes _kA$ is a
Cohen-Macaulay
ring provided $K$ is an algebraic extension field of $k$
and $A$ is a Cohen-Macaulay ring such that $K\otimes _kA$ is Noetherian.
Indeed, let $P\in $ Spec$(K\otimes _kA)$ and $p=P\cap A$. Since $K\otimes _kA$
is a flat integral extension of $A$,
$ht(P)=ht(p)$ \cite[Lemma 2.1]{w}. By Theorem 1.1, $G(p)=G(K\otimes
_kp)\leq G(P)$. Therefore $ht(P)=ht(p)=G(p)\leq G(P)\leq ht(P)$
\cite[Theorem 138]{k}.
Thus $G(P)=ht(P)$. Consequently, $K\otimes _kA$ is a Cohen-Macaulay ring.

Step 2. Let $K$ and $L$ be any extension fields of $k$ such that $
K\otimes _kL$ is Noetherian. We may suppose that $t=t.d.(K:k)<\infty $.
Let $x_1,...,x_t$ be elements of $K$ algebraically
independent over $k$. Then $K\otimes _kL \cong K\otimes
_{k(x_1,...,x_t)}S^{-1}L[x_1,...,x_t]$ (\cite[Proposition 2.6]{sh}), where $
S=k[x_1,...,x_t]\setminus \{0\}$. Since $K$ is an algebraic extension field of $
k(x_1,...,x_t)$ and $A=S^{-1}L[x_1,...,x_t]$ is a Cohen-Macaulay ring
(\cite[Theorem 151 and Theorem 139]{k}), by Step1, $K\otimes _kL \cong
K\otimes _{k(x_1,...,x_t)}A$
is a Cohen-Macaulay ring. $\Box $\bigskip

\noindent {\bf Proposition 2.3.} Let $A$ and $B$ be $k$-algebras such that
$A\otimes _kB$ is Noetherian.\ Let $P$ be a prime ideal of $A\otimes _kB$,
$p=P\cap A$, and $q=$ $P\cap B$. Then

a) $ht(P)=ht(p)+ht(q)+ht(\frac P{p\otimes_kB+A\otimes_kq})$.

b) $G(P(A\otimes _kB)_P)=G(pA_p)+G(qB_q)+ht(\frac P{p\otimes
_kB+A\otimes_kq})$.\bigskip

\noindent {\it Proof.} a) Consider the canonical flat homomorphism of
Noetherian rings \[f : A\rightarrow A\otimes_kB.\]
Applying \cite[Theorem 15.1]{m}, we have
\[ \begin{array}{ll}
ht(P)&=ht(p)+dim(\frac {(A\otimes_kB)_P}{p(A\otimes_kB)_P})\\
     &=ht(p)+dim(\frac {(A\otimes_kB)_P}{(p\otimes_kB)_P})\\
     &=ht(p)+dim((\frac {A\otimes_kB}{p\otimes_kB})_{\frac P{p\otimes_kB}})\\
     &=ht(p)+ht(\frac P{p\otimes_kB}).
\end{array}\]
Similarly, via the canonical homomorphism of Noetherian rings \[g:B\rightarrow \frac Ap\otimes_kB,\] we get
\[
\begin{array}{ll}
ht(\frac P{p\otimes_kB})&=ht(q)+ht(\frac {P/(p\otimes_kB)}{\frac Ap\otimes_kq})\\
&=ht(q)+ht(\frac {P/(p\otimes_kB)}{(p\otimes_kB+A\otimes_kq)/(p\otimes_kB)})\\
&=h(q)+ht(\frac P{p\otimes_kB+A\otimes_kq}).
\end{array}\]
It follows that $ht(P)=ht(p)+ht(q)+ht(\frac P{p\otimes_kB+A\otimes_kq})$, as desired.

b)  Notice first that $k_A(p)\otimes _kk_B(q)$ is a
Cohen-Macaulay ring, by Lemma 2.2.
Set $S_1 = A\setminus p$, $S_2 = B\setminus q$, and $S=\{a\otimes b|a\in S_1$ and $b\in
S_2\}$.

The above homomorphism $f$ induces the local flat homomorphism of Noetherian rings $A_p\rightarrow (A\otimes_kB)_P$.
In view of \cite[Corollary,
p. 181]{m}
or \cite[Proposition IV-6.3.1]{gr}, we have
\[ \begin{array}{ll}
G(P(A\otimes _kB)_P)  & = $ depth$ ((A\otimes _kB)_P) \\
                      & = $ depth$ (A_p)+$ depth$ (\frac{(A\otimes
_kB)_P}{pA_p(A\otimes _kB)_P}) \\
  & = G(pA_p)+\;$depth$ (\frac{(A\otimes _kB)_P}{(p\otimes _kB)_P}) \\
  & = G(pA_p)+$ depth$ ((\frac Ap\otimes_kB)_{\frac P{p\otimes _kB}}).
\end{array} \]

In a similar way, via the induced local flat homomorphism $B_q\rightarrow
(\frac Ap\otimes _kB)_{\frac P{p\otimes _kB}}$ of $g$, we get
\begin{center}
depth$ ((\frac Ap\otimes _kB)_{\frac P{p\otimes _kB}})=
G(qB_q)+$ depth$ ((\frac Ap\otimes
_k\frac Bq)_{\frac P{p\otimes _kB+A\otimes _kq}})$.
\end{center}
It follows that
\[ \begin{array}{ll}
G(P(A\otimes_kB)_P) & = G(pA_p)+ G(qB_q)+$ depth$ ((\frac Ap\otimes _k\frac Bq)_{\frac P{p\otimes_kB+A\otimes _kq}}) \\
                    & = G(pA_p)+G(qB_q)+$ depth$ ((k_A(p)\otimes _kk_B(q))_H) \\
                    & = G(pA_p)+G(qB_q)+\dim ((k_A(p)\otimes _kk_B(q))_H),
\end{array} \]
where $H = {\frac{S^{-1}P}{S_1^{-1}p\otimes _kS_2^{-1}B+S_1^{-1}A \otimes
_kS_2^{-1}q}}$.
Consequently, $G(P(A\otimes _kB)_P)=G(pA_p)+G(qB_q)+ht(\frac
P{p\otimes_kB+A\otimes _kq})$,
as desired. $\Box $\bigskip

\noindent {\it Proof of the theorem.} i)$\Rightarrow$ ii) and
ii)$\Rightarrow$ iii) are
obvious. Assume that (iii) holds. Let $p\in $ Spec$(A)$ and $q\in $ Spec$(B)$.
Then, by Theorem 1.1, $G(p\otimes_kB + A\otimes_kq) = G(p)+G(q)$. On the
other
hand, by Proposition 2.3, $ht(p\otimes_kB + A\otimes_kq) =ht(p)+ht(q)$.
Hence, since $G(p\otimes_kB + A\otimes_kq) =ht(p\otimes_kB +
A\otimes_kq)$, $G(p)+G(q)=ht(p)+ht(q)$. Therefore $ht(p)-G(p)=G(q)-ht(q)$,
so that $G(p)=ht(p)$ and $G(q)=ht(q)$, making (iv) hold.

Now, suppose that (iv) holds. By \cite[Theorem 140]{k}, it is sufficient to
prove that
$(A\otimes_kB)_P$ is a Cohen-Macaulay ring for each prime ideal $P$ of
$A\otimes_kB$.
Let $P$ be a prime ideal of $A\otimes_kB$, $p=P\cap A$, and $q=P\cap B$.
By Proposition 2.3, $G(P(A\otimes_kB)_P)
=G(pA_p)+G(qB_q)+ht(\frac{P}{p\otimes_kB+A\otimes_kq})$ and
$ht(P(A\otimes_kB)_P)=ht(P)=ht(p)+ht(q)+ht(\frac{P}{p\otimes_kB+A\otimes_kq})$.
Since $A$ and $B$ are Cohen-Macaulay, $A_p$ and $B_q$ are
Cohen-Macaulay. Then $G(pA_p)=ht(p)$ and $G(qB_q)=ht(q)$. Therefore
$G(P(A\otimes_kB)_P)= ht(P(A\otimes_kB)_P)$. Then
$(A\otimes_kB)_P$ is a Cohen-Macaulay ring. Hence (i) holds. $\Box$
\bigskip

\noindent {\bf Remark 2.4.} One may prove directly (i)$\Leftrightarrow$(iv) of Theorem 2.1 by using Lemma 2.2 and \cite[Corollaire IV-6.3.3]{gr}. However,
our proof is designed to draw extra benefits: Proposition 2.3 and Theorem 2.1(ii) \& (iii) shed more light on the prime ideal structure of (Noetherian) tensor products of $k$-algebras. Further, in the Noetherian case, Proposition 2.3(a) stands for a satisfactory analogue of [4, Theorem 1], a central result for polynomial
 rings. \bigskip

\noindent {\bf Remark 2.5.} Theorem 2.1 may
allow one to
determine the grade for new categories of primes of $A \otimes_kB$
(different from those treated in Theorem 1.1). For instance, let
$P\in $ Spec$(A\otimes _kB)$ with $p=P\cap A$
and $q=P\cap B$. Assume that $p$ and $q$ are generated by an
$A$-sequence and a $B$-sequence, respectively. Then
\[G(P)=G(p)+G(q)+G(\frac P{p\otimes _kB+A\otimes _kq}).\]
If, in addition, $p$ and $q$ are maximal ideals, then
\[G(P)=G(p)+G(q)+ht(\frac P{p\otimes_kB+A\otimes _kq}).\]

Indeed, let $p=(x_1,...,x_n)$ and $q=(y_1,...,y_m)$ such that
$x_1,...,x_n$ is an $A$-sequence and $y_1,...,y_m$ is a $B$-sequence.
Clearly, $x_1,...,x_n,y_1,...,y_m$ is an $A\otimes _kB$-sequence in
$p\otimes _kB+A\otimes _kq$ with $p\otimes _kB+A\otimes_kq
=(x_1,...,x_n,y_1,...,y_m)$. By \cite[Theorem 116]{k},
$G(\frac P{p\otimes _kB+A\otimes _kq})=G(P)-(n+m)=G(P)-G(p)-G(q)$.

Assume now that $p$ and $q$ are maximal ideals. Then, applying Theorem 2.1,
$\frac{A\otimes _kB}{p\otimes _kB+A\otimes _kq}\cong \frac Ap \otimes
_k\frac Bq$
is a Cohen-Macaulay ring. It follows that
$G(\frac P{p\otimes _kB+A\otimes _kq})=ht(\frac P{p\otimes _kB+A\otimes
_kq})$, as desired. \bigskip\bigskip


\end{document}